\documentclass{amsart}
\theoremstyle{plain}
\newtheorem{Thm}{Theorem}

\newtheorem{Prop}[Thm]{Proposition}
\newtheorem{Lem}[Thm]{Lemma}

\begin{document}
\large
\title[Radial Symmetry and Monotonicity Results for an Integral Equation]
{Radial Symmetry and Monotonicity Results for an Integral
Equation}

\author{Li MA and DeZhong Chen}

\address{Department of mathematical sciences \\
Tsinghua university \\
Beijing 100084 \\
China}

\email{lma@math.tsinghua.edu.cn} \dedicatory{}
\date{Aug 3rd, 2004}

\begin{abstract}
In this paper, we consider radial symmetry property of positive
solutions of  an integral equation arising from some higher order
semi-linear elliptic equations on the whole space $\mathbf{R}^n$.
We do not use the usual way to get symmetric result by using
moving plane method. The nice thing in our argument is that we
only need a Hardy-Littlewood-Sobolev type inequality. Our main
result is Theorem 1 below.
\end{abstract}

\thanks{The work is partially supported by 973 key project of science and technology ministry of
China. The first named author would like to thank Prof. W.X.Chen
for a helpful discussion in the summer of 2004 in Tsinghua
University, Beijing}
 \maketitle

\begin{section}{Introduction}

  In the study of standing waves of the non-linear Klein-Gordon
  equations we meet the following semi-linear elliptic equation
  (EE):
$$
  -\Delta u+u=u|u|^{\beta-1}, \quad{on}\quad{\mathbf{R}^n}
$$
where $\beta>1$ is a positive constant. This equation also appears
from the ground states of the Schrodinger equation \cite{B}. It is
shown that smooth positive solutions of (EE) are unique and radial
symmetric, so it has a nice decay at infinity \cite{Kw}. In this
paper, we consider radial symmetry property of positive solutions
of  an integral equation arising from some higher order
semi-linear elliptic equations on the whole space $\mathbf{R}^n$.
Our work is motivated from the works of Chen, Li, and Ou
\cite{CLO03} \cite{CLO04}, where the authors studied the Yamabe
type equations. The usual way to get symmetric result is using the
classical moving plane method. However, the moving plane method is
limited since it uses the Maximum principle. The nice thing
observed by Chen-Li-Ou \cite{CLO03} is that one only needs a
Hardy-Littlewood-Sobolev type inequality to get symmetric result
for positive solutions of elliptic equations of Yamabe type.

  For $\alpha>0$, $\beta>1$, we consider positive solutions for
  the following semi-linear partial differential equation in
  $\mathbf{R}^{n}$:
  \begin{equation}
  (I-\Delta)^{\alpha/2}(u)=u^{\beta}
  \end{equation}
  where $\Delta=\sum_{i=1}^{n}\partial^{2}/\partial x_{i}^{2}$
  denotes the Laplace operator in $\mathbf{R}^{n}$.
  With some decay assumption of solutions at infinity,
  it is well-known that $(1)$ is equivalent to the following
  integral equation
  \begin{equation}
  u=g_{\alpha}\ast u^{\beta},
  \end{equation}
  where $\ast$ denotes the convolution, and $g_{\alpha}$ is the Bessel kernel
  (for its precise definition, see $\S2$).

  Our main result is
  \begin{Thm}
  Assume that  $u\in L^{q}(\mathbf{R}^{n})$,
  where $q>\max\{\beta,n(\beta-1)/\alpha\}$ is a positive constant, is a positive
  solution of $(2)$. Then it must be radially symmetric and
  monotone decreasing about some point.
  \end{Thm}
\end{section}

\begin{section}{Bessel Potentials}
  For the convenience of readers, we recall some basic properties of Bessel
  potentials.

  The Bessel kernel, $g_{\alpha}$, $\alpha>0$, is defined by
  $$
  g_{\alpha}(x)=\frac{1}{\gamma(\alpha)}
  \int_{0}^{\infty}\exp(-\pi|x|^{2}/\delta)\exp(-\delta/4\pi)
  \delta^{(-n+\alpha)/2}\frac{d\delta}{\delta},
  $$
  where $\gamma(\alpha)=(4\pi)^{\alpha/2}\Gamma(\alpha/2)$.

  We now state, without proof, certain elementary facts about
  $g_{\alpha}$. For details, one may refer to ~\cite{S70} or
  ~\cite{Z89}.
  \begin{Prop}
  $ $

  (1) The Fourier transform of $g_{\alpha}$ is
  $$
  \widehat{g_{\alpha}}(x)=\frac{1}{(1+4\pi^{2}|x|^{2})^{\alpha/2}}
  $$
  where the Fourier transform is
  $$
  \widehat{f}(x)=\int_{\mathbf{R}^{n}}f(t)\exp(-2\pi ix\cdot t)dt.
  $$

  (2) For each $\alpha>0$, $g_{\alpha}(x)\in
  L^{1}(\mathbf{R}^{n})$.

  (3) The following Bessel composition formula holds
  $$
  g_{\alpha}\ast g_{\beta}=g_{\alpha+\beta},\alpha,\beta\geq0.
  $$
  \end{Prop}

  For any $\alpha\geq0$, and $f\in L^{p}(\mathbf{R}^{n})$, $1\leq
  p\leq\infty$, we define the Bessel potentials $B_{\alpha}(f)$ as
  $$
  B_{\alpha}(f)=g_{\alpha}\ast f
  $$
  if $\alpha>0$, and
  $$
  B_{0}(f)=f.
  $$

  \begin{Prop}
  $ $

  (1) $\|B_{\alpha}(f)\|_{L^{p}(\mathbf{R}^{n})}\leq\|f\|_{L^{p}(\mathbf{R}^{n})}$, $1\leq
  p\leq\infty$.

  (2) $B_{\alpha}\cdot
  B_{\beta}=B_{\alpha+\beta},\alpha,\beta\geq0$.
  \end{Prop}
  \begin{proof}
  $ $

  $(1)$ follows from the Young's inequality and the fact that
  $\|g_{\alpha}\|_{L^{1}(\mathbf{R}^{n})}=1(=\int_{\mathbf{R}^{n}}g_{\alpha}(x)dx)$.

  $(2)$ follows from the Bessel composition formula.
  \end{proof}
  We remark that on Sobolev spaces \cite{S70},
  $$
  B_{\alpha}=(I-\triangle)^{-\alpha/2}.
  $$

  The most interesting fact concerning Bessel potentials is that
  they can be employed to characterize the Sobolev spaces
  $W^{k,p}(\mathbf{R}^{n})$. This is expressed in the following
  theorem where we employ the notation
  $$
  L^{\alpha,p}(\mathbf{R}^{n}),\alpha>0,1\leq p\leq\infty,
  $$
  to denote all functions $u$ such that
  $$
  u=g_{\alpha}\ast f
  $$
  for some $f\in L^{p}(\mathbf{R}^{n})$. We define
  $$
  \|u\|_{L^{\alpha,p}(\mathbf{R}^{n})}=\|f\|_{L^{p}(\mathbf{R}^{n})}.
  $$
  With respect to this norm, $L^{\alpha,p}(\mathbf{R}^{n})$ is a
  Banach space.
  \begin{Thm}
  If $k$ is a nonnegative integer and $1<p<\infty$, then
  $$
  L^{k,p}(\mathbf{R}^{n})=W^{k,p}(\mathbf{R}^{n}).
  $$
  Moreover, if $u\in L^{k,p}(\mathbf{R}^{n})$ with $u=g_{\alpha}\ast f$,
  then
  $$
  C^{-1}\|f\|_{L^{p}(\mathbf{R}^{n})}\leq\|u\|_{W^{k,p}(\mathbf{R}^{n})}\leq
  C\|f\|_{L^{p}(\mathbf{R}^{n})}
  $$
  where $C=C(\alpha,n,p)$.
  \end{Thm}
\end{section}
\begin{section}{Imbedding Theorems}
  In this section, we prove some imbedding results which will play
  an important role in the proof of Theorem 1. First of all, let
  us recall the classical Sobolev imbedding theorem of the
  spaces $W^{k,p}(\mathbf{R}^{n})$. For the proof, see, for
  instance, ~\cite{A75}.
  \begin{Thm}(The Sobolev imbedding theorem) For $1\leq p<\infty$, there exist the
  following imbeddings:
  $$
  W^{k,p}(\mathbf{R}^{n}) \longrightarrow
  \{
    \begin{array}{lll}
    L^{q}(\mathbf{R}^{n}),\hspace{4mm} kp<n,\hspace{2mm} and\hspace{2mm} p\leq q\leq\frac{np}{n-kp},\\
    L^{q}(\mathbf{R}^{n}),\hspace{4mm} kp=n,\hspace{2mm} and\hspace{2mm} p\leq q<\infty,\\
    C_{b}(\mathbf{R}^{n}),\hspace{4mm} kp>n,\\
    \end{array}
  $$
  where $C_{b}(\mathbf{R}^{n})=\{u\in C(\mathbf{R}^{n}): u\hspace{1mm}is\hspace{1mm}bounded
  \hspace{1mm}on\hspace{1mm}\mathbf{R}^{n}\}$.
  \end{Thm}
  The following lemma comes from ~\cite{A75}(Theorem 7.63 (d) and
  (e)).
  \begin{Lem}
  $ $

  $(1).$ If $t\leq s$ and $1<p\leq
  q\leq\frac{np}{n-(s-t)p}<\infty$, then
  $$
  L^{s,p}(\mathbf{R}^{n})\rightarrow L^{t,q}(\mathbf{R}^{n}).
  $$

  $(2).$ If $0\leq\mu\leq s-n/p<1$, then
  $$
  L^{s,p}(\mathbf{R}^{n})\rightarrow C^{0,\mu}(\mathbf{R}^{n}).
  $$
  \end{Lem}
  Now we can prove
  \begin{Lem}
  Assume that $q>\max\{\beta,n(\beta-1)/\alpha\}$. If $f\in
  L^{q/\beta}(\mathbf{R}^{n})$, then $B_{\alpha}(f)\in
  L^{q}(\mathbf{R}^{n})$. Moreover, we have
  $$
  \|B_{\alpha}(f)\|_{L^{q}(\mathbf{R}^{n})}\leq C\|f\|_{L^{q/\beta}(\mathbf{R}^{n})}
  $$
  where $C=C(\alpha,\beta,n,q)$.
  \end{Lem}
  \begin{proof}
  Obviously, we have that $B_{\alpha}(f)\in L^{\alpha,q/\beta}(\mathbf{R}^{n})$
  since $f\in L^{q/\beta}(\mathbf{R}^{n})$. When $\alpha$ is an integer,
  by Theorem 4, $B_{\alpha}(f)\in
  W^{\alpha,q/\beta}(\mathbf{R}^{n})$. Then we can
  directly use Theorem 5 to obtain the results. We leave the proof
  as an exercise to the interested readers. Here, we only give the
  proof of the case when $\alpha$ is a fraction, i.e.,
  $\alpha-[\alpha]>0$, where $[\alpha]$ denotes the integer
  satisfying $[\alpha]\leq\alpha<[\alpha]+1$. There are two cases.

  {\bf Case 1.} $q<\frac{n\beta}{\alpha-[\alpha]}$

  Then by Lemma 6 (1), we have
  $$
  B_{\alpha}(f)\in L^{[\alpha],r}(\mathbf{R}^{n}),
  $$
  where $1<q/\beta\leq
  r\leq r_{0}<\infty$ and
  $r_{0}=\frac{nq/\beta}{n-(\alpha-[\alpha])q/\beta}$. Since
  $[\alpha]$ is an integer, we have
  $$
  B_{\alpha}(f)\in W^{[\alpha],r}(\mathbf{R}^{n}),
  $$

  There are two subcases.

  {\bf Case 1.1.} $q\leq r_{0}$

  This is equivalent to
  $$
  q\geq\frac{n(\beta-1)}{\alpha-[\alpha]}.
  $$

  In this case,  we have that $B_{\alpha}(f)\in W^{[\alpha],q}(\mathbf{R}^{n})$.
  Moreover,
  $$
  \|B_{\alpha}(f)\|_{L^{q}(\mathbf{R}^{n})}
  \leq C\|B_{\alpha}(f)\|_{L^{q/\beta}(\mathbf{R}^{n})}
  =C\|f\|_{L^{q/\beta}(\mathbf{R}^{n})}.
  $$

  {\bf Case 1.2.} $q>r_{0}$

  This is equivalent to
  $$
  q<\frac{n(\beta-1)}{\alpha-[\alpha]}.
  $$
  In this case, we need to use Theorem 5 to raise the exponent.

  When $q<n\beta/\alpha$, we have
  $$
  [\alpha]r_{0}<n.
  $$
  To make sure that $q$ is in
  $[r_{0},\frac{nr_{0}}{n-[\alpha]r_{0}}]$, it requires that
  $$
  q\geq n(\beta-1)/\alpha.
  $$
  Then by Theorem 5, we get
  $$
  B_{\alpha}(f)\in L^{q}(\mathbf{R}^{n}).
  $$

  When $q=n\beta/\alpha$, we have
  $$
  [\alpha]r_{0}=n.
  $$
  Then by Theorem 5, we get
  $$
  B_{\alpha}(f)\in L^{q}(\mathbf{R}^{n}).
  $$

  When $q>n\beta/\alpha$, we have
  $$
  [\alpha]r_{0}>n.
  $$
  Then by Theorem 5, we get
  $$
  B_{\alpha}(f)\in C_{b}(\mathbf{R}^{n}).
  $$
  A straightforward calculation shows that
  $$
  \|B_{\alpha}(f)\|_{L^{q}(\mathbf{R}^{n})}
  \leq\|B_{\alpha}(f)\|_{C_{b}(\mathbf{R}^{n})}^{1-\frac{1}{\beta}}
  \|B_{\alpha}(f)\|_{L^{q/\beta}(\mathbf{R}^{n})}^{\frac{1}{\beta}}
  $$
  $$
  \hspace{31mm}\leq C\|B_{\alpha}(f)\|_{L^{q/\beta}(\mathbf{R}^{n})}^{1-\frac{1}{\beta}}
  \|B_{\alpha}(f)\|_{L^{q/\beta}(\mathbf{R}^{n})}^{\frac{1}{\beta}}
  $$
  $$
  \hspace{3mm}=C\|B_{\alpha}(f)\|_{L^{q/\beta}(\mathbf{R}^{n})}
  $$
  $$
  \hspace{-4mm}=C\|f\|_{L^{q/\beta}(\mathbf{R}^{n})}.
  $$

  {\bf Case 2.} $q\geq\frac{n\beta}{\alpha-[\alpha]}$

  First, we see that
  $$
  B_{[\alpha]}(f)\in L^{\alpha-[\alpha],q/\beta}(\mathbf{R}^{n}).
  $$

  Notice that
  $$
  0\leq\alpha-[\alpha]-\frac{n}{q/\beta}<1
  \Leftrightarrow
  q\geq\frac{n\beta}{\alpha-[\alpha]}.
  $$
  So by Lemma 6 (2), we have
  $$
  B_{[\alpha]}(f)\in C^{0,\mu}(\mathbf{R}^{n}),
  $$
  where $0\leq\mu\leq\alpha-[\alpha]-\frac{n}{q/\beta}$. Then as
  before, we compute
  $$
  \|B_{[\alpha]}(f)\|_{L^{q}(\mathbf{R}^{n})}
  \leq\|B_{[\alpha]}(f)\|_{C^{0,\mu}(\mathbf{R}^{n})}^{1-\frac{1}{\beta}}
  \|B_{[\alpha]}(f)\|_{L^{q/\beta}(\mathbf{R}^{n})}^{\frac{1}{\beta}}
  $$
  $$
  \hspace{30mm}\leq C\|B_{[\alpha]}(f)\|_{L^{q/\beta}(\mathbf{R}^{n})}^{1-\frac{1}{\beta}}
  \|B_{[\alpha]}(f)\|_{L^{q/\beta}(\mathbf{R}^{n})}^{\frac{1}{\beta}}
  $$
  $$
  =C\|B_{[\alpha]}(f)\|_{L^{q/\beta}(\mathbf{R}^{n})}
  $$
  $$
  \hspace{-9mm}=C\|f\|_{L^{q/\beta}(\mathbf{R}^{n})}.
  $$

  By now, we have finished the proof.
\end{proof}
\end{section}
\begin{section}{Proof of Theorem 1}
  For a given real number $\lambda$, define
  $$
  \Sigma_{\lambda}=\{x=(x_{1},\cdots,x_{n})|x_{1}\geq\lambda\}.
  $$
  Let $x^{\lambda}=(2\lambda-x_{1},\cdots,x_{n})$,
  $u_{\lambda}(x)=u(x^{\lambda})$.
  \begin{Lem}
  For any solution $u(x)$ of $(2)$, we have
  $$
  u(x)-u_{\lambda}(x)=\int_{\Sigma_{\lambda}}(g_{\alpha}(x-y)-g_{\alpha}(x^{\lambda}-y))
  (u(y)^{\beta}-u_{\lambda}(y)^{\beta})dy.
  $$
  \end{Lem}
  \begin{proof}
  Let
  $$
  \Sigma_{\lambda}^{c}=\{x=(x_{1},\cdots,x_{n})|x_{1}<\lambda\}.
  $$
  Then it is easy to see that
  $$
  u(x)=\int_{\Sigma_{\lambda}}g_{\alpha}(x-y)u(y)^{\beta}dy
  +\int_{\Sigma_{\lambda}^{c}}g_{\alpha}(x-y)u(y)^{\beta}dy
  $$
  $$
  \hspace{13mm}=\int_{\Sigma_{\lambda}}g_{\alpha}(x-y)u(y)^{\beta}dy
  +\int_{\Sigma_{\lambda}}g_{\alpha}(x^{\lambda}-y)u_{\lambda}(y)^{\beta}dy.
  $$
  Substituting $x$ by $x^{\lambda}$, we get
  $$
  u(x^{\lambda})=\int_{\Sigma_{\lambda}}g_{\alpha}(x^{\lambda}-y)u(y)^{\beta}dy
  +\int_{\Sigma_{\lambda}}g_{\alpha}(x-y)u_{\lambda}(y)^{\beta}dy.
  $$
  Thus
  $$
  \hspace{-88mm}u(x)-u(x^{\lambda})
  $$
  $$
  =\int_{\Sigma_{\lambda}}g_{\alpha}(x-y)(u(y)^{\beta}-u_{\lambda}(y)^{\beta})dy
  -\int_{\Sigma_{\lambda}}g_{\alpha}(x^{\lambda}-y)(u(y)^{\beta}-u_{\lambda}(y)^{\beta})dy
  $$
  $$
  \hspace{-36mm}=\int_{\Sigma_{\lambda}}(g_{\alpha}(x-y)-g_{\alpha}(x^{\lambda}-y))(u(y)^{\beta}-u_{\lambda}(y)^{\beta})dy.
  $$
  \end{proof}
  {\bf Proof of Theorem 1.} Our proof is divided into two steps.

  {\bf Step 1.} Define
  $$
  \Sigma_{\lambda}^{-}=\{x|x\in\Sigma_{\lambda},u(x)<u_{\lambda}(x)\}.
  $$
  We want to show that for sufficiently negative values of
  $\lambda$, $\Sigma_{\lambda}^{-}$ must be empty.

  Whenever $x,y\in\Sigma_{\lambda}$, we have that $|x-y|\leq|x^{\lambda}-y|$. Then by the definition of
  $g_{\alpha}$, we have
  $$
  g_{\alpha}(x-y)\geq g_{\alpha}(x^{\lambda}-y).
  $$
  Then by Lemma 8, for $x\in\Sigma_{\lambda}^{-}$,
  $$
  u_{\lambda}(x)-u(x)\leq\int_{\Sigma_{\lambda}^{-}}(g_{\alpha}(x-y)-g_{\alpha}(x^{\lambda}-y))
  (u_{\lambda}(y)^{\beta}-u(y)^{\beta})dy
  $$
  $$
  \hspace{2mm}\leq\beta\int_{\Sigma_{\lambda}^{-}}g_{\alpha}(x-y)[u_{\lambda}^{\beta-1}(u_{\lambda}-u)](y)dy.
  $$
  It follows first from Lemma 7 and then the
  $H\ddot{o}lder$
  inequality that
  $$
  \|u_{\lambda}-u\|_{L^{q}(\Sigma_{\lambda}^{-})}\leq\beta\|B_{\alpha}(u_{\lambda}^{\beta-1}
  (u_{\lambda}-u))\|_{L^{q}(\Sigma_{\lambda}^{-})}
  $$
  $$
  \hspace{21mm}\leq C\|u_{\lambda}^{\beta-1}(u_{\lambda}-u)\|_{L^{q/\beta}(\Sigma_{\lambda}^{-})}
  $$
  $$
  \hspace{38mm}\leq C(\int_{\Sigma_{\lambda}^{-}}u_{\lambda}(y)^{q}dy)^{\frac{\beta-1}{q}}
  \|u_{\lambda}-u\|_{L^{q}(\Sigma_{\lambda}^{-})}
  $$
  \begin{equation}
  \hspace{35mm}\leq C(\int_{\Sigma_{\lambda}^{c}}u(y)^{q}dy)^{\frac{\beta-1}{q}}
  \|u_{\lambda}-u\|_{L^{q}(\Sigma_{\lambda}^{-})}
  \end{equation}
  Since $u\in L^{q}(\mathbf{R}^{n})$, we can choose $N$
  sufficiently large, such that for $\lambda\leq-N$, we have
  $$
  C(\int_{\Sigma_{\lambda}^{c}}u(y)^{q}dy)^{\frac{\beta-1}{q}}\leq\frac{1}{2}.
  $$
  Now $(3)$ implies that
  $$
  \|u_{\lambda}-u\|_{L^{q}(\Sigma_{\lambda}^{-})}=0,
  $$
  and therefore $\Sigma_{\lambda}^{-}$ must be measure zero, and
  hence empty.

  {\bf Step 2.} Now we have that for $\lambda\leq-N$,
  \begin{equation}
  u(x)\geq u_{\lambda}(x),\hspace{5mm}\forall
  x\in\Sigma_{\lambda}.
  \end{equation}
  Thus we can start moving the plane continuously from
  $\lambda\leq-N$ to the right as long as $(4)$ holds. Suppose
  that at a $\lambda_{0}<0$, we have $u(x)\geq
  u_{\lambda_{0}}(x)$, but $meas\hspace{0.5mm}\{x\in
  \Sigma_{\lambda_{0}}|u(x)>u_{\lambda_{0}}(x)\}>0$. We will show
  that the plane can be moved further to the right, i.e., there
  exists an $\epsilon$ depending on $n$, $\alpha$, $\beta$, $q$
  and the solution $u$ such that $u(x)\geq u_{\lambda}(x)$ on
  $\Sigma_{\lambda}$ for all $\lambda$ in
  $[\lambda_{0},\lambda_{0}+\epsilon)$.

  By Lemma 8, we see that $u(x)>u_{\lambda_{0}}(x)$ in the
  interior of $\Sigma_{\lambda_{0}}$. Let
  $$
  \overline{\Sigma_{\lambda_{0}}^{-}}=\{x\in
  \Sigma_{\lambda_{0}}|u(x)\leq u_{\lambda_{0}}(x)\}.
  $$
  Then it is easy to see that
  $\overline{\Sigma_{\lambda_{0}}^{-}}$ has measure zero, and
  $\lim_{\lambda\rightarrow\lambda_{0}}\Sigma_{\lambda}^{-}\subset\overline{\Sigma_{\lambda_{0}}^{-}}$.
  Let $(\Sigma_{\lambda}^{-})^{\star}$ be the reflection of
  $\Sigma_{\lambda}^{-}$ about the plane $x_{1}=\lambda$. From the
  third inequality of $(3)$, we get
  \begin{equation}
  \|u_{\lambda}-u\|_{L^{q}(\Sigma_{\lambda}^{-})}\leq C
  (\int_{(\Sigma_{\lambda}^{-})^{\star}}u(y)^{q}dy)^{\frac{\beta-1}{q}}
  \|u_{\lambda}-u\|_{L^{q}(\Sigma_{\lambda}^{-})}.
  \end{equation}
  Since $u\in L^{q}(\mathbf{R}^{n})$, then one can choose
  $\epsilon$ small enough, such that for all $\lambda$ in
  $[\lambda_{0},\lambda_{0}+\epsilon)$,
  $$
  C(\int_{(\Sigma_{\lambda}^{-})^{\star}}u(y)^{q}dy)^{\frac{\beta-1}{q}}\leq\frac{1}{2}.
  $$
  Now by $(5)$, we have
  $$
  \|u_{\lambda}-u\|_{L^{q}(\Sigma_{\lambda}^{-})}=0,
  $$
  and therefore $\Sigma_{\lambda}^{-}$ must be empty.
\end{section}

\end{document}